\documentclass[twoside]{book}
\usepackage{graphicx}
\usepackage{epsf}
\usepackage{ilcmacro,url}

 \makeatletter
\ifx\UNDEF\mail\def\mail{ }\else\fi
\ifx\UNDEF\prange\def\prange{0 0}\else\fi

\gdef\@empty{}
\def\Mail#1 #2 {\gdef\thecontact{#1}\gdef\theaddr{#2}}
\def\Range#1 #2 {\gdef\thefirstpage{#1}\gdef\thelastpage{#2}}
{\let\'\mail \expandafter\Mail\' }  
{\let\'\prange \expandafter\Range\' }   
 \gdef\@shtitle{\relax}
 \long\def\shtitle#1{\gdef\@shtitle{#1}}
 \long\def\author#1{\gdef\@author{#1}}
 
 \gdef\@abstract{}
 \long\def\abstract#1{\gdef\@abstract{#1}}

 \def\maketitle{\thispagestyle{empty}\chapter{\@title}}
 \renewcommand\chapter{\if@openright\cleardoublepage\else\clearpage\fi
                    \thispagestyle{empty}%
                    \global\@topnum\z@
                    \@afterindentfalse
                    \secdef\@chapter\@schapter}
 \def\@makechapterhead#1{%
  \vspace*{50\p@}%
  {\parindent \z@ \raggedleft \normalfont
    \ifnum \c@secnumdepth >\m@ne
      \if@mainmatter
        \par\nobreak
        \vskip 20\p@
      \fi
    \fi
    \interlinepenalty\@M
    \Huge \bfseries #1\par\nobreak
    \vskip.25in
    \large\bfseries\@author\par\nobreak
    \vskip 40\p@}
    \ifx\@abstract\@empty\else{\small\@abstract\par\vskip20\p@}\fi
  }

\DeclareRobustCommand\em
        {\@nomath\em \ifdim \fontdimen\@ne\font >\z@
                       \upshape \else \slshape \fi}

\def\@begintheorem#1#2{\sl \trivlist \item[\hskip \labelsep{\bf #1\ #2}]}
\def\@opargbegintheorem#1#2#3{\sl \trivlist
     \item[\hskip \labelsep{\bf #1\ #2\ (#3)}]}

 \newcommand{\sectlabel}[1]{\label{sect:#1}}

 \renewcommand{\theequation} {\mbox{\arabic{equation}}}
 \renewcommand{\thesection} {\mbox{\arabic{section}}}

  \def\@arabic#1{\number #1} 

\long\def\@makecaption#1#2{
    \vskip\abovecaptionskip
    \sbox\@tempboxa{{\small {\bf #1}: #2}}%
    \ifdim\wd\@tempboxa>\hsize
        {\small {\bf #1}: #2\par}
    \else
       \global\@minipagefalse
       \hbox to\hsize{\hfil\box\@tempboxa\hfil}
    \fi
    \vskip \belowcaptionskip}

\def\figstrut#1{\hbox to\linewidth{\vrule height#1\hfill}}

\renewenvironment{thebibliography}[1]
     {\section*{\bibname
        \@mkboth{\MakeUppercase\bibname}{\MakeUppercase\bibname}}%
      \list{\@biblabel{\@arabic\c@enumiv}}%
           {\settowidth\labelwidth{\@biblabel{#1}}%
            \leftmargin\labelwidth
            \advance\leftmargin\labelsep
            \@openbib@code
            \usecounter{enumiv}%
            \let\p@enumiv\@empty
            \renewcommand\theenumiv{\@arabic\c@enumiv}}%
      \sloppy
      \clubpenalty4000
      \@clubpenalty \clubpenalty
      \widowpenalty4000%
      \sfcode`\.\@m}
     {\def\@noitemerr
       {\@latex@warning{Empty `thebibliography' environment}}%
      \endlist}
\makeatother

\shtitle{Random graph models of communication network topologies}
\title{Random graph models of communication network topologies}
\author{Hannu Reittu and Ilkka Norros\\ first.last@vtt.fi}
\date{21 September 2007}
\abstract{We consider a variant of so called power-law random graph. A sequence of expected degrees corresponds to a power-law degree distribution with finite mean and infinite variance. In previous works the asymptotic picture with number of nodes limiting to infinity has been considered. It was found that an interesting structure appears. It has resemblance with such graphs like the Internet graph. Some simulations have shown that a finite sized variant has similar properties as well. Here we investigate this case in more analytical fashion, and, with help of some simple lower bounds for large valued expectations of relevant random variables, we can shed some light into this issue.  A new term, 'communication range random graph' is introduced to emphasize that some further restrictions are needed to have a relevant random graph model for a reasonable sized communication network, like the Internet. In this case a pleasant model is obtained, giving the opportunity to understand such networks on an intuitive level. This would be beneficial in order to understand, say, how a particular routing works in such networks.  }
\begin{document}
\maketitle
\section{Introduction}\sectlabel{intro}
Since the pioneering work by three Faloutsos brothers, \cite{3Faloutsos99} and some other groups like that of Barabasi's \cite{barabasialbert02}, around the millennium, so called power-law graphs, with degree sequence obeying power-law distribution with finite mean and infinite variance, have attracted high interest by several authors. This degree sequence is argued to reflect some fundamental aspects of some communication networks and some other natural and technological networks. It appears that this concept can be turned into a mathematical object in several ways, with provable properties, see for instance \cite{bollobasjansonriordan07,aiellochunglu2000,chunglu03,hofhooznam-infvar,hooghiemstramieghem03,reittunorros06}. However, for the sake of tractability, asymptotic regime with growing number of nodes, $N\rightarrow\infty$, has been most popular. Some quite interesting results were obtained \cite{reittunorrosglobecom02,reittunorrosperfev04}. For instance, it was found that a random graph variant of this model, where with a given degree sequence links are drawn as randomly as possible, produces a graph with some characteristic properties that has correspondence with such network as the Internet itself, at its autonomous systems level (AS-graph). A kind of 'soft-hierarchy' of large degree nodes is formed only due to combinatorial probabilities. It is convenient to categorize nodes with increasing degrees into 'tiers' with nodes in certain subsequent intervals of degrees. It is sufficient to consider approximately $\log\log N$ tiers. This can be done in a way that a node in a tier has at least one link to upper-laying next tier, with probability tending to $1$, asymptotically. The hierarchical part of nodes was called the core. Thus very short ($\log\log N$) paths exist from the bottom to top of the core. It was also shown that even shorter paths are needed to find this core for almost any node in the same component. Thus, a $\log\log N$ asymptotical upper bound for distance in the giant component was established. These results with extensions were obtained independently by Chung and others, \cite{aiellochunglu2000,chunglu03}, using some refined methods of random graph theory. However, we found our more elementary approach with the concept of the core also very useful giving an intuitive insight. Refined variants of the theory can be found in works of van der Hofstad and others, \cite{hofhooznam-infvar}, showing, among other issues, that the $\log\log N$ upper bound is the best possible . 

However, it was apparent from the proofs, that convergence can be very slow, involving such functions like $1/\log\log N$ with limit $0$, that is approached only with 'unrealistically' large values of $N$, certainly unimaginable in the framework of communication networks. Some simulations indicated that in spite of this some reasonably sized graphs have properties that are similar to this asymptotic picture. That is why it is interesting to study this question in more details using an analytical approach. Here we do some first steps into this direction. 

It is also interesting, that the asymptotic model indicates interesting consequences for such graphs robustness against targeted attacks against the top level nodes. Such graphs show good robustness against such failures, at least in the terms of the distance: the remaining tiers are able to maintain connectivity with the price of only an insignificant number of extra hops. However, here it is also important to know how valid these results are for reasonable graph sizes. Recently, possibilities of extending the basic model by 'redirecting' the links, have been discussed. Here it is possible that the asymptotic picture is 'non convergent', meaning that it tells nothing about the finite variant of the graph.

\section{Model definition and asymptotic results}\sectlabel{model}

We consider a variant of power-law random graph, similar to one proposed by Chung and Lu \cite{chunglu03}, see also: \cite{bollobasjansonriordan07,reittunorros06} . A natural number, $N > 0$, is the number of nodes in the graph. Nodes are labeled with natural numbers $1,2,\cdots,N$. A node with label $i$, $1\leq i\leq N$ has 'capacity', $\lambda_i=(N/i)^\alpha$, with real number $\alpha$, ${1\over2}<\alpha< 1$, which reflects the power-law degree sequence. For each possible unordered pair of nodes $\{i,j\}$, $i,j\in \{1,2,\cdots,N\}$, we associate the number of links between those nodes as a random variable $E_{i,j}$ with Poisson distribution, with expected value $\E(E_{i,j})=\lambda_i \lambda_j/\sum_{i=1}^N \lambda_i$. All these random variables are considered as independent. In shorthand we write $E_{i,j}\cong Po(\lambda_i \lambda_j/L_N)$, with $L_N\equiv \sum_{i=1}^N \lambda_i$. Thus multiple links and self-loops are allowed. However, such artifacts are not too harmful, since the vast majority of these variables take values $0$ or $1$ in a large enough graph. The expected degree of node $i$, $d_i$, is thus $d_i= \E(\sum_{j=1}^N E_{i,j})=\sum_{j=1}^N \lambda_i\lambda_j/\sum_{i=1}^N \lambda_i=\lambda_i$, due to a basic property for the sum of independent Poisson distributed random variables. Thus $\lambda_i$ has the meaning of expected degree of node $i$.

Let us define the following sequence of functions:
\begin{eqnarray}
\beta_0(N)&=&{1\over \tau -1}+{\epsilon(N)\over \tau-2}\\
\nonumber
\beta_j(N)&=&(\tau-2)\beta_{j-1}(N)+\epsilon(N),\quad j=1,2,\cdots
\end{eqnarray}
with ${1\over \tau -1}=\alpha$, $\epsilon(N)=l(N)/\log N$ and $l(\cdot)$ is a very slowly diverging function as its argument grows to infinity. 

We define the 'upper layers' as 
$$
U_0\equiv \{1\},\quad U_j\equiv \{i: \lambda_i\geq N^{\beta_j(N)}\},\quad j=1,2,\cdots.
$$
Provided that $l(\cdot)$ fulfills: $l(1)=1$, $l(N)/\log\log\log N \rightarrow 0$ and $l(N)/\log\log\log\log N \rightarrow \infty$, we had the following result for the power-law graph described above:
\begin{theorem}\label{loglogN}
Let
$$
k^*\equiv k^*(N):=\left\lceil {\log\log N\over -\log(\tau-2)}\right\rceil,
$$
where $\left\lceil \cdot\right\rceil$, denotes the least integer greater than or equal to its argument.
Then the hop-count distance between two randomly chosen vertices of the giant component, which exists asymptotically almost surely (a.a.s.), is less than $2 k^*(N)(1+o(1))$, a.a.s.  
\end{theorem}
We define the core, $C$, as the upper layer $U_{k^*}$. Later on this proposition was strengthened considerably, one such result being that this upper bound is tight \cite{hofhooznam-infvar}. However, such more detailed analysis is very involved and that is why we prefer to stay on the level of simple upper bounds, also in what follows. Roughly speaking, the idea of our proof of Theorem \ref{loglogN} is that the probability that a node $i$ in any layer $U_j$ has a link to the upper layer $U_{j-1}$ with probability tending to $1$ as $N\rightarrow \infty$ (we write this as: $\pr{ U_j \ni i\leftrightarrow U_{j-1}}\rightarrow 1$). Thus, it takes at most $k^*$ hops to travel from the lowest layer $U_k^*$ to the top degree node. Further, almost any node within the giant component has a path to some node in $U_{k^*}$, with a number of hops that is sub-linear with $\log\log N$. Thus, the upper bound follows. However, as we see in the next section, the convergence of probability to $1$ can be very slow, say, $\pr{ U_j\ni i \leftrightarrow U_{j-1}}\geq 1-{c\over \log\log N}, \quad c>0$, a convergence rate that is practically 'unobservable' in our framework. Thus it is a relevant question whether this asymptotic picture tells anything about a graph with only a reasonably large number of nodes. Some simulations seem to indicate that the answer is positive, see e.g. \cite{tangmunarunkitwillinger02,reittunorrosperfev04, norrosreittu03}. However, we found that in order to find a corresponding above described layer structure in a finite model one must define function $l(\cdot)$ in some particular way, not prescribed by its asymptotic behavior only. In this paper our aim is to explain such circumstances and indicate a way how such finite sized random graphs can be analysed, and thus to make such random graph models more usable for modeling communication networks. 

\section{Analysis of 'communication range' graphs}\label{communicationgraph}

It is easy to see that the cardinality of layer $j>0$ is: $\mid U_j\mid= \left\lfloor N^{1-\beta_j(N)/\alpha}\right\rfloor= \left\lfloor N^{1-(\tau-1)\beta_j(N)}\right\rfloor$, where, $\left\lfloor \cdot \right\rfloor$ is the largest integer smaller or equal to the argument. Thus we have a lower bound for the sum of capacities in a layer $j>0$:
$$
V(U_j)\equiv\sum_{i\in U_j}\lambda_i\geq N^{\beta_j(N)}\left\lfloor   N^{1-(\tau-1)\beta_j(N)}\right\rfloor \geq N^{1-(\tau-2)\beta_j(N)}-N^{\beta_j(N)}\equiv V_0(U_j)
$$
and $V(U_0)=N^\alpha$. For $L_N$ we have asymptotically linear scaling with $N$. Indeed, it is easy to see that $L_N\geq N^\alpha\int_0^N(1+x)^{-\alpha}dx\geq {1\over 1-\alpha}(N-N^\alpha)\geq {c\over 1-\alpha}N$, where $c$ can be taken arbitrarily close to $1$, provided that $N$ is large enough. For instance, $c=9/10$ is valid provided that $N>10^{1\over 1-\alpha}$. However, notably these bounds are not uniform with $\alpha$. For $\alpha$, close to $1$, we would have to choose a small value for $c$, for any reasonable $N$. This is a general trend here, since we must also fix the range of $\alpha$ more precisely, not just stating that $1/2<\alpha <1$, which was sufficient for the asymptotic analysis. This circumstance reflects the fact that the asymptotic regime is approached sensitively with respect to $\alpha$. Our hypothesis is that for communication networks with $N$ in reasonable range of thousands of nodes or tens thousands of nodes, it is necessary to have $\alpha$ in the lower half of the range $(0,1)$  then the asymptotic range is reasonably close. Luckily, in the case of the Internet, this is a range of $\alpha$ that has been observed. We call this range of $N$ and $\alpha$, the \textit{communication range}. Similarly we have:
$$       
V_0(U_j)\geq c_j N^{1-(\tau-2)\beta_j(N)}
$$ 
with constants $c_j$, close to $1$, provided $N$ and $\tau$ are in the communication range. Within the same range, we can easily find a lower bound for the probability that a node in layer $j$ has at least one link to layer $j-1$:
\begin{eqnarray}\label{pj} 
\pr{ U_j(N)\ni i\leftrightarrow U_{j-1}(N)}=1-\exp(-\lambda_i V(U_{j-1})/L_N)\\\nonumber\geq 1-\exp(-{c_{j-1}\over c}e^{l(N)})\equiv p_j.
\end{eqnarray}  
 This relation also shows the delicacy of the communication range, where we can approximate $c_j/c$ by $1$ --- otherwise we would need a number far from $1$, giving a big effect to the lower bound $p_0$. Say, if $\alpha$ is close to $1$, this ratio would be a big number resulting in a very low probability. Notably, $p_0$ is also sensitive to the choice of function $l(\cdot)$. In the asymptotic sense these circumstances are irrelevant, since in any case $p_0\rightarrow 1$ as $N\rightarrow\infty$, if only $l(N)\rightarrow\infty$.

We wish to show that in communication range, a particularly defined 'core' has a similar role as in the asymptotic graph. In particular, the $\log\log N$ scaling of distance is roughly valid. We show that the lower bound of expected number of nodes that have a link to a core node that has a path through the core hierarchy up to the top is large enough and suggests that the $\log\log N$ upper bound is valid within the communication range, for the vast majority of nodes.

Assume fixed $N$ and take a natural number $x>0$. The probability that a node in layer $U_x$, has a link to $U_{x-1}$ is lower bounded by $p_x$. The probability that the same node has a path to $U_{x-2}$, through $U_{x-1}$, is lower bounded by $p_x p_{x-1}$. And so forth, probability that the same node has a path to the top node $1$, going through at most $x$ layers, is lower bounded by $p_x p_{x-1}\cdots p_1$. Denoting by $c_x$ the minimal ratio  $c_j\over c$ in relation (\ref{pj}), in the corresponding range of $j$, and denoting
\begin{eqnarray}\label{p0}
p_0= 1-\exp(-c_xe^{(3-\tau)l(N)}),
\end{eqnarray}   
we find that the probability that a node in layer $U_x$ has a path described above is lower bounded by $p_0^x$. Note that, in asymptotic range, $p_0$ tends to $1$ (very slowly), however, within our finite range this is an important parameter affecting the quality of bounds. 

Denote by $U_x'\subset U_x$ the nodes in layer $U_x$ having a valid path with upmost $x$ hops within the subsequent layers to the top node $1$. As a result we have:
\begin{eqnarray}\label{xlayervolume}
\E{\mid U_x'\mid}\geq p_0^x \mid U_x\mid\,\\\nonumber
\E{V(U_x')}\geq N^{\beta_x}p_0^x \mid U_x\mid.
\end{eqnarray}
Denote by $N_x$ the nodes that are in $U_x'$ or have a link to a node in it. The probability that a node outside $U_x'$ has link to it is lower bounded by $1-\exp(-{V(U_x')\over L_N})\geq {1\over 2}{V(U_x')\over L_N}$. That is why, for the conditional expectation, we have:
$$
\E(\mid\ N_x\mid\mid\mid U_x'\mid)\geq{N\over 2} {V(U_x')\over L_N}.
$$ 
Therefore, according to (\ref{xlayervolume}),
\begin{eqnarray}\label{xdensity}
\E(\mid N_x\mid)\geq {N^{\beta_x}N^1\over 2 L_N}p_0^x \mid U_x\mid.
\end{eqnarray}

The task is to maximise the bound (\ref{xdensity}), in a way that $x$ is not too large. It appears that asymptotically we end up with the setting that is in line with Theorem \ref{loglogN}. However, it is also possible to find a 'setting' of $l(\cdot)$ and $x$ that corresponds to asymptotic-like behavior in the communication range. To get a qualitative picture, we make simplifications in relation (\ref{xdensity}), assuming all constants, that are close to $1$, equal to one. As a result we get an approximate lower bound for $\E(\mid N_x\mid)\geq s(x,l)$:
\begin{eqnarray}\label{s(x,l)}
s(x,l)\approx N^{1-{(\tau-2)^{x+1}\over\tau-1}}(1-\exp(-m))^x m^{-{\tau-2\over 3-\tau}},\quad m\equiv\exp((3-\tau)l).
\end{eqnarray}
The maximum is found as solution of equations ${\partial s(x,l)\over \partial x}=0, {\partial s(x,l)\over \partial l}=0 $, and by comparing values of the function in the closest integer arguments. Although the equations are not solvable in closed form, the first one yields the relation:
\begin{eqnarray}\label{max}
{\partial s(x,l)\over \partial x}=0\Leftrightarrow x={\log\log N\over -\log(\tau-2)}+\log\left({\log(1/(1-\exp(-m))\over -\log(\tau-2)}\right)-1,
\end{eqnarray}
The first term is analogous to the one in Theorem \ref{loglogN}, the next one depends on the choice of $l$ through $m(l)$, which should be found from the second equation. However, we can see the asymptotical regime from these equations. Indeed, we see that the relevant factor of $s$, with respect to argument $l$, is asymptotically $m^{-{\tau-2\over 3-\tau}}\exp\left({-x\over \exp (m)}\right)\sim m^{-{\tau-2\over 3-\tau}}\exp\left({\log\log N\over \log(\tau-2)\exp (m)}\right)  $, where we took into account the equation (\ref{max}). The second term suggests that $m$ should be an increasing function, at least threefold logarithm, and $l$ should be no slower than fourfold logarithm. The first term suggests that it should not be too fast, in this respect the lowest possible would be the best, and gives the maximum. However, as suggested by the Theorem \ref{loglogN}, the leading term is indifferent with respect to this range. Indeed, if we make the corresponding substitution to (\ref{max}), we see that the leading term is just $k^*$.

\begin{figure}
\begin{minipage}{5.9cm}
\includegraphics[width=5.9cm]{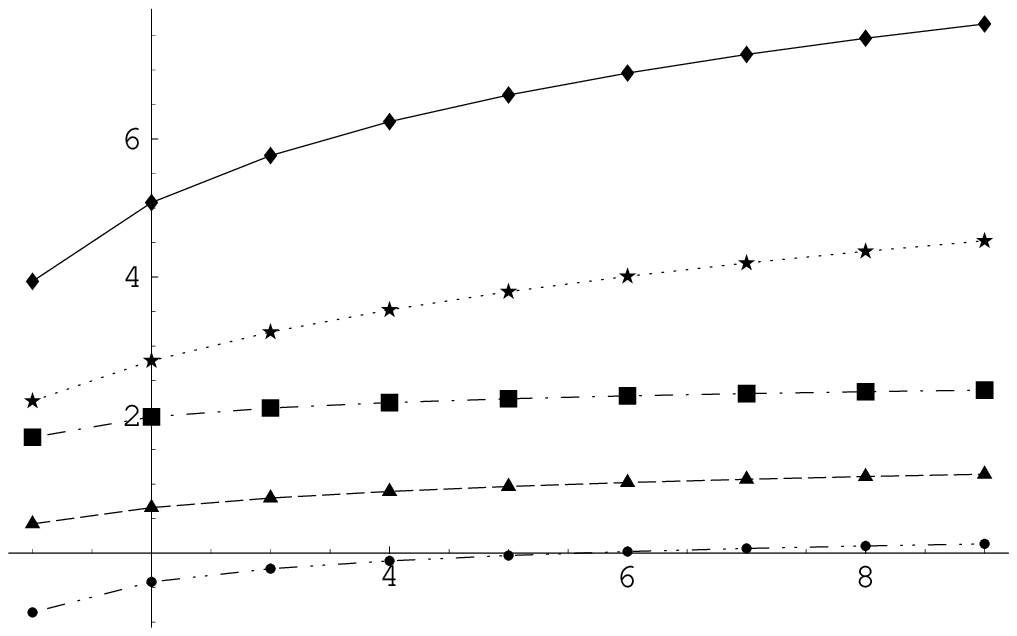}
\caption{Plot of functions, calculated in points $N=10^{k+1}, k=1,2,\cdots, 9$: $x(N)$, ${1\over \log{1/(\tau-2)}}\log\log{N} $, $l(N)$, $\log\log\log N$ and $\log\log\log\log N$, listed from top to down.}
\label{kuva1}
\end{minipage}
\hspace*{2mm}
\begin{minipage}{5.9cm}
\includegraphics[width=5.9cm]{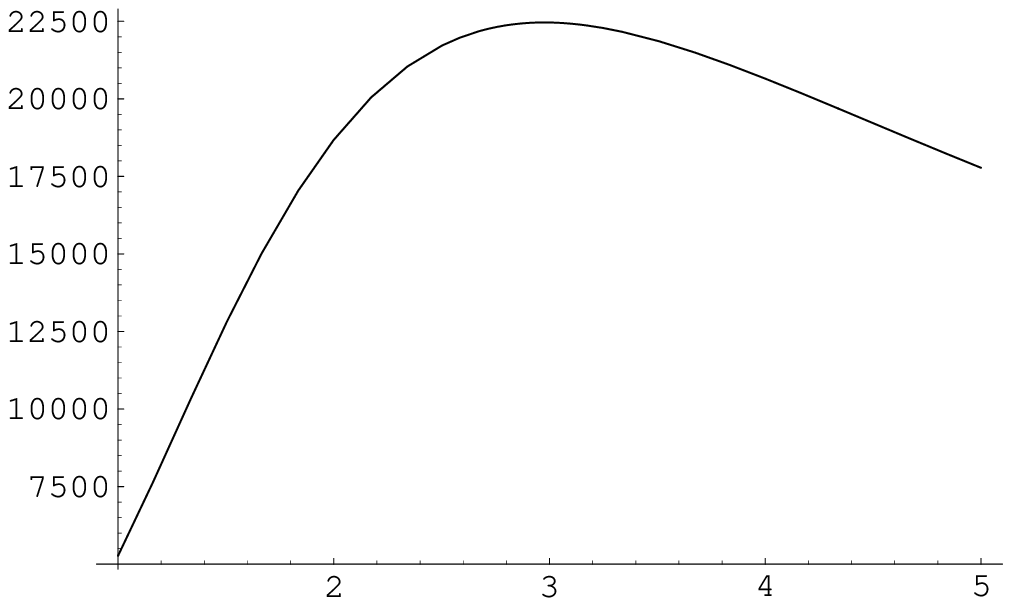}
\caption{Plot of function $s(x^*,l)$, with $N=10^5$ and with fixed first argument with value $x^*$, that it takes at maximum. }
\label{kuva2}
\end{minipage}
\end{figure}

If we substitute these asymptotic estimates, as arguments to (\ref{s(x,l)}), we find
\begin{eqnarray}\label{asydensity}
s(x,l)\propto N \exp(-(\tau-2)l(N)),
\end{eqnarray}
which is only slightly lower than $N$. Thus, the factor $\exp(-(\tau-2)l(N))$ has the meaning of lower bound expected 'density' of neighbors of the core nodes, with valid paths to the top node. This density is almost constant, as a function of $N$, between $1/\log\log\log N$ and $1/\log\log N$. In the communication range, one would guess that the best choice would be $1/\log\log\log N$. Numerical calculations seems to support this, see Figure  \ref{kuva1}. Our next plot, in Figure \ref{kuva2}, indicates that in this range the choice of this function has some effect. 
By taking for the $x$ its optimal value $x^*$, we see that the function $s(x^*,l)$, see \ref{s(x,l)}, has a visible maximum. This numerics showed also that in a very wide range of $N$, from $100$ to $10^{10}$, a rather constant fraction of nodes, in lower bound expectation, are neighbors of core with valid paths to the top. In our case with $\tau=2.5$, this was around $20$ percent. As a result, a random node is able to find a node in the core that has a path to its top, with a moderate sized search: approximately every fifth node is of this type.

This kind of 'quasi-stationarity' or extremely low dependence on $N$ should be  good news, since it simplifies the usability of such models. For instance, $s(x,l)$ and $l$ can be taken as constant parameters hardly changing in any reasonable range. We also see that the height of the core is almost constant, and its major term is a function of type $k^*(N)$. One drawback is that we have only lower bound type results and some unrigorous estimates were done. However, it is quite likely that a thorough analysis will not reveal any substantial new features, although it is mandatory to check it. It would be interesting to compare this approach with a 'conceptual model' of the Internet, called the 'Jellyfish', \cite{jellyfish06}. 

{\bf Acknowledgments} This work was financially supported by EU-project Net-Refound, Project Number: 034413.    
\bibliographystyle{apt}

\end{document}